# Contractive upper triangular matrices with prescribed diagonals and superdiagonals

Axel RENARD

**Abstract**

A characterization of the matrix representation of the compressed shift acting on a finite-dimensional model space endowed with the Takenaka-Malmquist-Walsh basis, among all upper triangular matrices, is proved. For a fixed dimension, this matrix is the unique matrix with spectral norm no greater than one and with a prescribed diagonal and superdiagonal. We also discuss an extension to the setting of infinite matrices.

## 1 Introduction and main results

For any $n \geq 2$ and given $n$ points $\omega_1, \ldots, \omega_n$ in the open unit disk $\mathbb{D}$, we consider the $n \times n$ matrix $M_n$ defined by:

$$[M_n]_{i,j} = \begin{cases} \omega_j & \text{if } i=j \\ \prod_{k=i+1}^{j-1} (-\overline{\omega_k})\sqrt{1-|\omega_i|^2}\sqrt{1-|\omega_j|^2} & \text{if } i<j \\ 0 & \text{if } i>j \end{cases}. \qquad (1.1)$$

For instance, for $n = 2$ and $n = 3$, we have:

$$M_2 = \begin{pmatrix} \omega_1 & \sqrt{1-|\omega_1|^2}\sqrt{1-|\omega_2|^2} \\ 0 & \omega_2 \end{pmatrix},$$

$$M_3 = \begin{pmatrix} \omega_1 & \sqrt{1-|\omega_1|^2}\sqrt{1-|\omega_2|^2} & -\overline{\omega_2}\sqrt{1-|\omega_1|^2}\sqrt{1-|\omega_3|^2} \\ 0 & \omega_2 & \sqrt{1-|\omega_2|^2}\sqrt{1-|\omega_3|^2} \\ 0 & 0 & \omega_3 \end{pmatrix}.$$

The matrix $M_n = M_n(\omega_1, \cdots, \omega_n)$ will be referred as the *model matrix of size $n$*, with eigenvalues $\omega_1, \ldots, \omega_n$. Indeed, $M_n$ is the matrix representation of the *compressed shift* acting on the $n$-dimensional model space when endowed with the Takenaka-Malmquist-Walsh basis. Detailed definitions will be provided in Section 2. Interested readers can find further information and a more comprehensive overview of the topic in [7, 8, 6].

Those matrices $M_n$ arise in various mathematical problems. It seems that their first occurrence was in the work of Ptak and Young [16, 21], where they appeared as extremal matrices for the problem of finding

$$C(H_n, p, f) = \sup\{\|f(A)\| : A \in \mathcal{B}(H_n), \|A\| \leq 1, p(A) = 0\},$$

where $H_n$ is a $n-$dimensional Hilbert space, $p$ is a monic polynomial of degree $n$ whose roots lie in the open unit disk, and $f$ is a function analytic in the neighborhood of the roots of $p$. The same matrices also appear in another extremal problem studied in [13].

Furthermore, those matrices $M_n$ are also related to matrices of class $S_n$. Recall (see [20, Chapter 7]) that a $n \times n$ matrix $A$ is said to be of class $S_n$ whenever $A$ is a contraction (*i.e.* of spectral norm no greater than one), with all eigenvalues in the open unit disk, and such that the rank of $I - A^*A$ is one. It is known [9] that $A$ is of class $S_n$ if and only if it is unitarily similar to a matrix model $M_n = M_n(\omega_1, \cdots, \omega_n)$ for some points $\omega_1, \ldots, \omega_n$ in $\mathbb{D}$. Some spectral properties of these matrices, and different properties of their numerical range, have been studied in [20, 10, 12, 2, 3] and the references therein.

Another recent appearance of the matrices $M_n$ was in [1], where criteria for matrix contractivity and their applications to Schwarz-Pick inequalities were discussed. For $2 \times 2$ matrices, a contractivity criterion is easily derived by noting that the spectral norm (also called Euclidean operator norm) $\|T\|$ corresponds to the largest singular value of $T$. In particular, a straightforward computation shows that, given $\omega_1, \omega_2 \in \mathbb{D}$, the matrix $T_2 = \begin{pmatrix} \omega_1 & \alpha \\ 0 & \omega_2 \end{pmatrix} \in \mathcal{M}_2(\mathbb{C})$ is a contraction if and only if $|\alpha|^2 \leq (1 - |\omega_1|^2)(1 - |\omega_2|^2)$. In this sense, the matrix $M_2(\omega_1, \omega_2)$ can be regarded as an extremal case for this criterion, with the previous inequality becoming an equality.

A criterion to determine whenever a $3 \times 3$ upper-triangular matrix is a contraction has been established in [1]. This criterion leads, in particular, to the following observation: given $\omega_1, \omega_2, \omega_3 \in \mathbb{D}$, the matrix $M_3(\omega_1, \omega_2, \omega_3)$ can be regarded as an extremal case for $3 \times 3$ upper-triangular contractive matrices. More precisely, if we consider matrices of the form
$$T_3 = \begin{pmatrix} \omega_1 & \alpha_1 & \beta \\ 0 & \omega_2 & \alpha_2 \\ 0 & 0 & \omega_3 \end{pmatrix} \in \mathcal{M}_3(\mathbb{C}),$$
and if we fix $\alpha_1 = \sqrt{1 - |\omega_1|^2}\sqrt{1 - |\omega_2|^2}$, $\alpha_2 = \sqrt{1 - |\omega_2|^2}\sqrt{1 - |\omega_3|^2}$, then, $T_3$ is a contraction if and only if $T_3 = M_3 = M_3(\omega_1, \omega_2, \omega_3)$.

In this paper, we extend this observation, establishing the following theorem for $n \times n$ matrices:

**Theorem 1.1.** *Let $\omega_1, \ldots, \omega_n \in \mathbb{D}$ and let*

$$T_n = \begin{pmatrix} \omega_1 & \alpha_1^{(1)} & \alpha_1^{(2)} & \cdots & \cdots & \alpha_1^{(n-1)} \\ 0 & \omega_2 & \alpha_2^{(1)} & \alpha_2^{(2)} & & \vdots \\ \vdots & & \ddots & \ddots & \ddots & \vdots \\ & & & & & \alpha_{n-2}^{(2)} \\ & & & & & \alpha_{n-1}^{(1)} \\ 0 & \cdots & \cdots & \cdots & 0 & \omega_n \end{pmatrix} \in \mathcal{M}_n(\mathbb{C}).$$

*Assume that*
$$\alpha_i^{(1)} = \sqrt{1-|\omega_i|^2}\sqrt{1-|\omega_{i+1}|^2}, \text{ for all } 1 \leq i \leq n-1.$$

*Then, $T_n$ is a contraction if and only if $T_n = M_n$, where $M_n$ is defined as in* (1.1).

In other words, the matrix $T_n$ with prescribed diagonal and superdiagonal is a contraction if and only if

$$\alpha_i^{(j)} = \prod_{k=i+1}^{j+i-1}(-\overline{\omega_k})\sqrt{1-|w_i|^2}\sqrt{1-|\omega_{i+j}|^2}, \text{ for all } 1 \leq i \leq n, \ 1 \leq j \leq n-i.$$

Finally, we extend this result to infinite matrices, proving the following corollary:

**Corollary 1.2.** *Let $M = (m_{i,j})_{i,j=1}^{+\infty}$ be an infinite matrix of the form*

$$M = \begin{pmatrix} \omega_1 & \alpha_1^{(1)} & \alpha_1^{(2)} & \cdots \\ 0 & \omega_2 & \alpha_2^{(1)} & \alpha_2^{(2)} \\ 0 & & \ddots & \\ \vdots & & & \end{pmatrix},$$

*with:*

- $\omega_i \in \mathbb{D}$, *for all $i \geq 1$;*
- $\alpha_i^{(1)} = \sqrt{1-|\omega_i|^2}\sqrt{1-|\omega_{i+1}|^2}$, *for all $i \geq 1$.*

*Then, $M$ defines a contraction on $\ell^2(\mathbb{N})$ if and only if*

$$\alpha_i^{(j)} = \prod_{k=i+1}^{j+i-1}(-\overline{\omega_k})\sqrt{1-|w_i|^2}\sqrt{1-|\omega_{i+j}|^2}, \text{ for all } i,j \geq 1. \tag{1.2}$$

In particular, if the sequence $(\omega_k)_{k\geq 1}$ satisfies the Blaschke condition $\sum_{k=1}^{+\infty}(1-|\omega_k|) < +\infty$ (see Section 2), the last results means that $M$ defines a contraction on $\ell^2(\mathbb{N})$ if and only if it is the matrix representation of a compressed shift acting on an infinite dimensional model space. This can be related to the Sz.-Nagy-Foiaş model theory (see [7, Section 9.4] or [17]): if $T$ is a contraction acting on a Hilbert space satisfying $\|T^n x\| \to 0$, for all $x \in H$, and such that $\text{rank}(\text{Id} - T^*T) = \text{rank}(\text{Id} - TT^*) = 1$, then, $T$ is unitarily equivalent to a compressed shift.

## 2 Preliminaries about model space operator matrices

### 2.1 Notation

We start by introducing some notation. In this paper, we will denote by $\mathbb{D}$ the open unit disk, and by $\mathbb{T}$ the unit circle, *i.e.* $\mathbb{T} = \overline{\mathbb{D}}\backslash\mathbb{D}$. For $n \in \mathbb{N}^*$, we will also denote by $\mathcal{M}_n(\mathbb{C})$ the set of $n \times n$ matrices with complex coefficients, and by $\|\cdot\|$ the Euclidean operator norm on $\mathbb{C}^n$. Throughout this paper, $\mathbb{C}^n = (\mathbb{C}^n, \langle \cdot, \cdot \rangle)$ is always considered as an $n$-dimensional Hilbert space, where $\langle \cdot, \cdot \rangle$ denotes the standard inner product on $\mathbb{C}^n$. If $H$ and $K$ are two Hilbert spaces, we denote by $\mathcal{B}(H, K)$ the set of bounded linear operators from $H$ to $K$. We usually write $\mathcal{B}(H)$ instead of $\mathcal{B}(H, H)$. Moreover, if an operator $T$ satisfies $\|T\| \leq 1$, we say that $T$ is a *contraction*. Then, $T$ is a contraction if and only if $\mathrm{Id} - T^*T$ is positive semi-definite. Here Id is the identity operator and $T^*$ is the conjugate transpose of the matrix $T$, *i.e.* the (Hilbertian) adjoint of $T$. When $T$ is a contraction, we will denote by $D_T$ its defect operator defined by $D_T = (\mathrm{Id} - T^*T)^{1/2}$, where the *square root* denotes the unique positive semi-definite square root. Finally, in this manuscript, $[\![N, M]\!]$ denotes all the integers from $N$ to $M$, and we work with the usual convention that an empty sum is equal to 0, and an that empty product is equal to 1.

### 2.2 A review of model spaces

Let $H^\infty = H^\infty(\mathbb{D})$ be the set of all holomorphic functions that are bounded on $\mathbb{D}$. For $f \in H^\infty$, Fatou's theorem implies that $f$ has a well-defined radial boundary value $f(\zeta)$, for almost every $\zeta \in \mathbb{T}$. A function $u \in H^\infty$ is said to be *inner* if $|u(\zeta)| = 1$ almost everywhere on $\mathbb{T}$.

Now, let $H^2 = H^2(\mathbb{D})$ be the Hardy-Hilbert space of the unit disk, *i.e.* the set of analytic functions on $\mathbb{D}$ whose Taylor coefficients are square-summable, and let $\langle \cdot, \cdot \rangle_{H^2}$ denote the canonical scalar product on $H^2$. If $u$ is an inner function, the corresponding *model space* $\mathcal{K}_u$ is defined to be

$$\mathcal{K}_u := H^2(\mathbb{D}) \ominus uH^2(\mathbb{D}) = \left(uH^2\right)^\perp = \{f \in H^2(\mathbb{D}) \,:\, \langle f, uh \rangle_{H^2} = 0, \, \forall\, h \in H^2(\mathbb{D})\}.$$

Let $S$ be the unilateral shift $S: H^2 \ni f \mapsto zf \in H^2$. We define the associated *compressed shift* by $S_u := P_u S|_{\mathcal{K}_u}$, where $P_u$ is the orthogonal projection from $H^2(\mathbb{D})$ onto $\mathcal{K}_u$.

Now, let $\Theta_n = \prod_{k=1}^n \frac{z - \omega_k}{1 - \overline{\omega_k} z}$ be a Blaschke product with zeros $\omega_1, \ldots, \omega_n \in \mathbb{D}$ (repeated according to multiplicity), and let $b_{\omega_k}(z) = \frac{z - \omega_k}{1 - \overline{\omega_k} z}$ denote a single Blaschke factor. Define the functions $\phi_1, \ldots, \phi_n$ by

$$\phi_1(z) = \frac{\sqrt{1 - |\omega_1|^2}}{1 - \overline{\omega_1} z} \quad \text{and} \quad \phi_k(z) = \left(\prod_{j=1}^{k-1} b_{\omega_j}\right) \frac{\sqrt{1 - |\omega_k|^2}}{1 - \overline{\omega_k} z}, \quad k = 2, \ldots, n.$$

Then $(\phi_1, \ldots, \phi_n)$ is an orthonormal basis of $\mathcal{K}_{\Theta_n}$, called the Takenaka-Malmquist-Walsh basis of $\mathcal{K}_{\Theta_n}$.

Writing $S_{\Theta_n}$ with respect to the Takenaka-Malmquist-Walsh basis gives the matrix rep-resentation $M_n = M_n(\omega_1, \ldots, \omega_n)$, where:

$$[M_n]_{i,j} = \langle S_{\Theta_n} \phi_j, \phi_i \rangle$$
$$= \begin{cases} \omega_j & \text{if } i=j \\ \prod_{k=i+1}^{j-1} (-\overline{\omega_k}) \sqrt{1-|\omega_i|^2}\sqrt{1-|\omega_j|^2} & \text{if } i<j \\ 0 & \text{if } i>j \end{cases}.$$

Note that $M_n$ is a contraction, as it is the matrix representation of a compressed shift in an orthonormal basis.

This can be generalized to infinite Blaschke products. Let $\Theta_\infty$ be an infinite Blaschke product with zeros $\omega_1, \omega_2, \cdots \in \mathbb{D}$ repeated according to multiplicity. Recall that its zeros must satisfy the Blaschke condition $\sum_{k \geq 1}(1 - |\omega_k|) < \infty$. Then, the Takenaka-Malmquist-Walsh basis $(\phi_1, \phi_2, \ldots)$, defined by:

$$\phi_1(z) = \frac{\sqrt{1-|\omega_1|^2}}{1-\overline{\omega_1}z} \quad \text{and} \quad \phi_k(z) = \left(\prod_{j=1}^{k-1} b_{\omega_j}\right) \frac{\sqrt{1-|\omega_k|^2}}{1-\overline{\omega_k}z}, \quad k \geq 2,$$

is an orthonormal basis of $\mathcal{K}_{\Theta_\infty}$. The matrix representation $M_\infty = M_\infty((\omega_k)_{k\geq 1})$ of the compressed shift $S_{\Theta_\infty}$ in this basis is given by:

$$[M_\infty]_{i,j} = \langle S_{\Theta_\infty} \phi_j, \phi_i \rangle_{H^2}$$
$$= \begin{cases} \omega_j & \text{if } i=j \\ \prod_{k=i+1}^{j-1} (-\overline{\omega_k}) \sqrt{1-|\omega_i|^2}\sqrt{1-|\omega_j|^2} & \text{if } i<j \\ 0 & \text{if } i>j \end{cases}.$$

## 3 Extremal contractive matrices

In this section, we prove Theorem 1.1. The key ingredient of the proof will be the following result about matrix completion which goes back to Parrott [15] (see also [22, Theorem 12.22], [5], or [1, Appendix] for a revisited proof):

**Theorem 3.1** (Parrott). *Let $H_1, H_2, K_1, K_2$ be Hilbert spaces, and suppose that the operators $\begin{bmatrix} A \\ C \end{bmatrix} \in \mathcal{B}(H_1, K_1 \oplus K_2)$ and $\begin{bmatrix} C & D \end{bmatrix} \in \mathcal{B}(H_1 \oplus H_2, K_2)$ are contractions.*

*Then, $T = \begin{bmatrix} A & B \\ C & D \end{bmatrix} : H_1 \oplus H_2 \to K_1 \oplus K_2$ is a contraction if and only if there exists a contraction $W \in \mathcal{B}(H_2, K_1)$ such that*

$$B = D_{Z^*} W D_Y - Z C^* Y,$$

*where $Z \in \mathcal{B}(H_1, K_1)$ and $Y \in \mathcal{B}(H_2, K_2)$ are contractions such that $D = D_{C^*} Y$ and $A = Z D_C$.*

We will also need to use Möbius transformations in order to simplify the computations. For $\omega \in \mathbb{D}$ and $z \neq \frac{1}{\overline{\omega}}$, we define the Möbius transformation $M_\omega$ by

$$M_\omega(z) = \frac{\omega - z}{1 - \overline{\omega}z}.$$

We recall ([14, Chapter IX, Section 2]) that $M_\omega$ is involutive (that is, $M_\omega(M_\omega(z)) = z$), holomorphic on $\mathbb{D}$, and that for all $z \in \overline{\mathbb{D}}$, $|M_\omega(z)| \leq 1$, with equality if and only if $|z| = 1$. Using the rational functional calculus, we can define $M_\omega(T)$ when $\frac{1}{\overline{\omega}} \notin \sigma(T)$. In particular, this is the case for any contraction.

**Lemma 3.2.** *Let $\omega \in \mathbb{D}$. Suppose that the spectrum of the linear operator $T$ does not contain $\frac{1}{\overline{\omega}}$. Then*

(i) *The operator $M_\omega(M_\omega(T))$ is well defined and $M_\omega(M_\omega(T)) = T$;*

(ii) *The operator $T$ is a contraction if and only if $M_\omega(T)$ is a contraction.*

*Proof.* For the first part, we note that the spectrum $\sigma(M_\omega(T))$ of $M_\omega(T)$, which is equal to $M_\omega(\sigma(T))$ by the spectral mapping theorem, does not contain $\frac{1}{\overline{\omega}}$. This is because the equation $\frac{\omega - z}{1 - \overline{\omega}z} = \frac{1}{\overline{\omega}}$ has no solution. Therefore, $M_\omega(M_\omega(T))$ is well-defined. Subsequently, using properties of the rational functional calculus, we can write $M_\omega(M_\omega(T)) = (M_\omega \circ M_\omega)(T) = T$.

Using the previous item, for the second part it is enough to prove that if $T$ is a contraction, then so is $M_\omega(T)$. This can be seen as a consequence of the von Neumann inequality (see [19, 18]) or from the following elementary computation.
Let $x \in H$, and let $y = (\text{Id} - \overline{\omega}T)^{-1} x$. Then

$$||x||^2 - ||M_\omega(T)x||^2 = ||(\text{Id} - \overline{\omega}T)y||^2 - ||(\omega\text{Id} - T)y||^2$$
$$= \left(1 - |\omega|^2\right)\left(||y||^2 - ||Ty||^2\right) \geq 0.$$

Thus, $M_\omega(T)$ is a contraction. □

Now, we are ready for the proof of Theorem 1.1:

*Proof of Theorem 1.1.* First of all, note that if $T_n = (t_{i,j})_{1 \leq i,j \leq n}$, then we have $\alpha_i^{(j)} = t_{i,i+j}$. We will proceed by induction on the size of the matrix $T_n$.

- For $n = 2$, there is nothing to prove.

- Let $n \geq 3$, and assume that the result is true for a matrix of size $n - 1$.

If $T_n$ is a contraction, then, the two compressions

$$S_n^{(1)} = \begin{bmatrix} \omega_1 & \alpha_1^{(1)} & \alpha_1^{(2)} & \cdots & \alpha_1^{(n-2)} \\ 0 & \omega_2 & \alpha_2^{(1)} & \cdots & \alpha_2^{(n-3)} \\ \vdots & \ddots & \ddots & \ddots & \vdots \\ \vdots & & \ddots & \ddots & \alpha_{n-2}^{(1)} \\ 0 & \cdots & \cdots & 0 & \omega_{n-1} \end{bmatrix}$$

and

$$S_n^{(2)} = \begin{bmatrix} \omega_2 & \alpha_2^{(1)} & \alpha_2^{(2)} & \cdots & \alpha_2^{(n-2)} \\ 0 & \omega_3 & \alpha_3^{(1)} & \cdots & \alpha_3^{(n-3)} \\ \vdots & \ddots & \ddots & \ddots & \vdots \\ \vdots & & \ddots & \ddots & \alpha_{n-1}^{(1)} \\ 0 & \cdots & \cdots & 0 & \omega_n \end{bmatrix}$$

are contractions too. Using the induction hypothesis, that happens if and only if

$$\alpha_i^{(j)} = \prod_{k=i+1}^{j+i-1} (-\overline{\omega_k}) \sqrt{1 - |w_i|^2} \sqrt{1 - |\omega_{i+j}|^2}, \text{ for all } i, j,$$

except maybe for $\alpha_1^{(n-1)}$ (who does not appear in any of those two compressions). Our aim is now to prove that there exists a unique choice for the coefficient $\alpha_1^{(n-1)}$ such that $T_n$ is a contraction. As we know that for

$$a_1^{(n-1)} = \prod_{k=2}^{n-1} (-\overline{\omega_k}) \sqrt{1 - |w_1|^2} \sqrt{1 - |\omega_n|^2},$$

$T_n = M_{\Theta_n}$ is a contraction, that will be enough to prove the theorem.

Let $A = \begin{bmatrix} \omega_1 & \alpha_1^{(1)} & \alpha_1^{(2)} & \cdots & \alpha_1^{(n-2)} \end{bmatrix}$, $B = \begin{bmatrix} \alpha_1^{(n-1)} \end{bmatrix}$, $C = \begin{bmatrix} 0 & \omega_2 & \alpha_2^{(1)} & \cdots & \alpha_2^{(n-3)} \\ \vdots & \ddots & \ddots & \ddots & \vdots \\ \vdots & & \ddots & \ddots & \alpha_{n-2}^{(1)} \\ \vdots & & & \ddots & \omega_{n-1} \\ 0 & \cdots & \cdots & \cdots & 0 \end{bmatrix}$

and $D = \begin{bmatrix} \alpha_2^{(n-2)} \\ \vdots \\ \alpha_{n-1}^{(1)} \\ \omega_n \end{bmatrix}$. If $T_n$ is a contraction, then, $\begin{bmatrix} A \\ C \end{bmatrix}$ and $\begin{bmatrix} C & D \end{bmatrix}$ are contractions too. As it is a necessary condition for $T_n$ to be a contraction, we assume in the following that this condition is satisfied. By Parrott's theorem, $T_n$ is a contraction if and only if there exists a contraction $V$ such that:

$$B = (\text{Id} - ZZ^*)^{1/2} V (\text{Id} - Y^*Y)^{1/2} - ZC^*Y, \tag{3.1}$$

where $Y$ and $Z$ are contractions such that $D = (\mathrm{Id} - CC^*)^{1/2}Y$ and $A = Z(\mathrm{Id} - C^*C)^{1/2}$. Note that (3.1) can be rewritten:

$$\left|\alpha_1^{(n-1)} + ZC^*Y\right|^2 \leq (1 - ZZ^*) \times (1 - Y^*Y). \tag{3.2}$$

*Fact* 1. Suppose $\omega_2 = 0$. Then,

$$C^*C = \begin{bmatrix} 0 & 0 \cdots\cdots\cdots 0 \\ 0 & 0 & \ddots & \vdots \\ \vdots & \ddots & 1 & \ddots \\ \vdots & & \ddots & 0 \\ 0 \cdots\cdots\cdots & 0 & 1 \end{bmatrix} \begin{matrix} \updownarrow 2 \\ \\ \updownarrow n-3 \end{matrix}$$
$$\underbrace{\phantom{xx}}_{2} \underbrace{\phantom{xxxxx}}_{n-3}$$

*Proof.* Denote $C = (c_{i,j})_{1 \leq i,j \leq n-1}$, with $c_{i,j} = \begin{cases} 0 & \text{if } j \leq i \\ \omega_j & \text{if } j = i+1 \\ \alpha_{i+1}^{(j-i-1)} & \text{if } j \geq i+2 \end{cases}$.

We have $[C^*C]_{i,j} = \sum_{k=1}^{n-1} \overline{c_{k,i}} c_{k,j} = \sum_{k=1}^{m(i,j)} \overline{c_{k,i}} c_{k,j}$, where $m(i,j) = \min\{i-1, j-1\}$.

We can easily notice that $[C^*C]_{1,j} = 0$, for all $j \in [\![1, n-1]\!]$ and $[C^*C]_{i,1} = 0$, for all $i \in [\![1, n-1]\!]$ (which means that the first row and the first column of $C^*C$ are zero).

Moreover, it is easy to see that $[C^*C]_{2,2} = |\omega_2|^2 = 0$ and that $[C^*C]_{2,j} = \overline{\omega_2}\alpha_2^{(j-2)} = 0$, for all $j \in [\![3, n-1]\!]$.

Furthermore, for $i \in [\![3, n-1]\!]$, we get:

$$[C^*C]_{i,i} = \sum_{k=1}^{i-2} \left|\alpha_{k+1}^{(i-k-1)}\right|^2 + |\omega_i|^2$$
$$= (1 - |\omega_i|^2)\left[\sum_{k=1}^{i-2}\left(\prod_{s=k+2}^{i-1}|\omega_s|^2\right)(1-|\omega_{k+1}|^2)\right] + |\omega_i|^2$$
$$= (1 - |\omega_i|^2)\left[\sum_{k=1}^{i-2}\prod_{s=k+2}^{i-1}|\omega_s|^2 - \prod_{s=k+1}^{i-1}|\omega_s|^2\right] + |\omega_i|^2.$$

We recognize a telescopic sum, and we get:

$$[C^*C]_{i,i} = (1 - |\omega_i|^2)\left(1 - \prod_{s=2}^{i-1}|\omega_s|^2\right) + |\omega_i|^2$$
$$= 1,$$

as the first factor of the product $\prod_{s=2}^{i-1}|\omega_s|^2$ is equal to 0.

Similarly, for $i \in [\![3, n-1]\!]$ and $j > i$, we get:

$$[C^*C]_{i,j} = \sum_{k=1}^{i-2} \overline{\alpha_{k+1}^{(i-k-1)}} \alpha_{k+1}^{(j-k-1)} + \overline{\omega_i}\alpha_i^{(j-i)}$$

$$= \sqrt{1-|\omega_i|^2}\sqrt{1-|\omega_j|^2}\left[\sum_{k=1}^{i-1}\left((1-|\omega_{k+1}|^2)\prod_{s=k+2}^{i-1}(-\omega_s)\prod_{s=k+2}^{j-1}(-\overline{\omega_s})\right) - \prod_{s=i}^{j-1}(-\overline{\omega_s})\right]$$

$$= \sqrt{1-|\omega_i|^2}\sqrt{1-|\omega_j|^2}\left[\sum_{k=1}^{i-1}\left((1-|\omega_{k+1}|^2)\prod_{s=k+2}^{i-1}|\omega_s|^2\prod_{s=i}^{j-1}(-\overline{\omega_s})\right) - \prod_{s=i}^{j-1}(-\overline{\omega_s})\right]$$

$$= \sqrt{1-|\omega_i|^2}\sqrt{1-|\omega_j|^2}\left[\prod_{s=i}^{j-1}(-\overline{\omega_s})\right]\left[\sum_{k=1}^{i-1}\left(\prod_{s=k+2}^{i-1}|\omega_s|^2 - \prod_{s=k+1}^{i-1}|\omega_s|^2\right) - 1\right]$$

$$= \sqrt{1-|\omega_i|^2}\sqrt{1-|\omega_j|^2}\left[\prod_{s=i}^{j-1}(-\overline{\omega_s})\right]\left[1 - \prod_{s=2}^{i-1}|\omega_s|^2 - 1\right]$$

$$= 0.$$

Finally, we conclude using the self-adjointness of $C^*C$. □

*Fact* 2. Suppose $\omega_2 = 0$. Then $||T_n|| \leq 1$ if and only if $\alpha_1^{(n-1)} = 0$.

*Proof.* If $\omega_2 = 0$, we have

$$\text{Id} - C^*C = \begin{bmatrix} 1 & 0 & \cdots & \cdots & 0 \\ 0 & 1 & \ddots & & \vdots \\ \vdots & \ddots & 0 & \ddots & \vdots \\ \vdots & & \ddots & \ddots & 0 \\ 0 & \cdots & \cdots & 0 & 0 \end{bmatrix} \begin{matrix} \uparrow \\ 2 \\ \downarrow \\ \uparrow \\ n-3 \\ \downarrow \end{matrix}$$

$$\underbrace{\qquad}_{2} \underbrace{\qquad}_{n-3}$$

Denote $Z = \begin{bmatrix} z_1 & \cdots & z_{n-1} \end{bmatrix}$. The identity $A = Z(\text{Id} - C^*C)^{1/2}$ is equivalent to

$$\begin{cases} z_1 = \omega_1 \\ z_2 = \sqrt{1-|\omega_1|^2} \end{cases}.$$

The only possibility for $Z$ to remain a contraction is to have

$$z_3 = \cdots = z_{n-1} = 0.$$

Finally, we obtain:
$$ZZ^* = 1$$
and
$$ZC^* = \begin{bmatrix} 0 & \cdots & 0 \end{bmatrix}.$$

Therefore (3.2) holds if and only if:
$$\alpha_1^{(n-1)} = -ZC^*Y = 0.$$

The proof of Fact 2 is completed. □

The general case can be reduced to the case where $\omega_2 = 0$ by using the fact that $T_n$ is a contraction if and only if $M(T_n)$ is a contraction, where $M = M_{\omega_2} : z \mapsto \frac{\omega_2 - z}{1 - \overline{\omega_2} z}$ (Lemma 3.2). Using the explicit rational functional calculus (see *e.g.* [4, 11]), we can write:

$$M(T_n) = \begin{pmatrix} M(\omega_1) & \beta_1^{(1)} & * & \cdots & \cdots & * \\ 0 & M(\omega_2) & \beta_2^{(1)} & & & \vdots \\ \vdots & \ddots & \ddots & \ddots & & \vdots \\ \vdots & & \ddots & \ddots & \ddots & * \\ \vdots & & & \ddots & \ddots & \beta_{n-1}^{(1)} \\ 0 & \cdots & \cdots & \cdots & 0 & M(\omega_n) \end{pmatrix},$$

where $\beta_i^{(1)} = \alpha_i^{(1)} \cdot \frac{M(\omega_i) - M(\omega_{i+1})}{\omega_i - \omega_{i+1}} = \frac{\alpha_i^{(1)}(|\omega_2|^2 - 1)}{(1 - \overline{\omega_2}\omega_i)(1 - \overline{\omega_2}\omega_{i+1})}$.

Assume that
$$\alpha_i^{(1)} = \sqrt{1 - |\omega_i|^2}\sqrt{1 - |\omega_{i+1}|^2}, \quad \text{for all } i \in [\![1, n-1]\!].$$

Then, for all $i \in [\![1, n-1]\!]$, we have:
$$\beta_i^{(1)} = e^{i\phi_i}\sqrt{1 - |M(\omega_i)|^2}\sqrt{1 - |M(\omega_{i+1})|^2}, \quad \text{for some } \phi_i \in [0, 2\pi[.$$

Set $\theta_i = \phi_1 + \cdots + \phi_i$ and $U = \text{diag}(e^{i\theta_1}, \ldots, e^{i\theta_{n-1}}, 1)$. Then $U$ is unitary and

$$U^*M(T_n)U = \begin{pmatrix} M(\omega_1) & \gamma_1^{(1)} & * & \cdots & \cdots & * \\ 0 & M(\omega_2) & \gamma_2^{(1)} & & & \vdots \\ \vdots & \ddots & \ddots & \ddots & & \vdots \\ \vdots & & \ddots & \ddots & \ddots & * \\ \vdots & & & \ddots & \ddots & \gamma_{n-1}^{(1)} \\ 0 & \cdots & \cdots & \cdots & 0 & M(\omega_n) \end{pmatrix},$$

where $\gamma_i^{(1)} = \sqrt{1 - |M(\omega_i)|^2}\sqrt{1 - |M(\omega_{i+1})|^2}$.

The operator $M(T_n)$ is a contraction if and only if $U^*M(T_n)U$ is a contraction and, from the previous reasoning we know that if $U^*T_n U$ is a contraction, then, all its coefficients are uniquely determined. If all the coefficients of $U^*T_n U$ are uniquely determined, then, it is also the case for $M(T_n)$ and for $T_n = M \circ M(T_n)$. □

Now, we prove Corollary 1.2, extending the previous result to infinite matrices:

*Proof of Corollary 1.2.* Assume first that $M$ defines a contraction on $\ell^2(\mathbb{N})$. Let $\Theta_\infty$ be a Blaschke product with zeros $\omega_1, \omega_2, \ldots$ (repeated according to multiplicity). As $\mathcal{K}_{\Theta_\infty}$ is a separable Hilbert space, it is unitarily isomorphic to $\ell^2(\mathbb{N})$. Thus, $M$ defines a contraction on $\mathcal{K}_{\Theta_\infty}$, equipped with the Takenaka-Mamlquist-Walsh basis $(\phi_k)_{k=1}^{+\infty}$. Observe that $\text{span}(\phi_1, \ldots, \phi_n)$ is the model space $\mathcal{K}_{\Theta_n}$ associated to the finite Blaschke product $\Theta_n$ with zeros $\omega_1, \ldots, \omega_n$ (repeated according to multiplicity). For $n \geq 2$, let $P_n$ be the orthogonal projection from $\mathcal{K}_{\Theta_\infty}$ onto $\mathcal{K}_{\Theta_n}$. Then, every compression $P_n M P_n$ is a contraction and, by Theorem 1.1, we must have, for all $n \geq 2$,

$$\alpha_i^{(j)} = \prod_{k=i+1}^{j+i-1} (-\overline{\omega_k})\sqrt{1-|w_i|^2}\sqrt{1-|\omega_{i+j}|^2}, \text{ for all } 1 \leq i \leq n,\ 1 \leq j \leq n-i,$$

which is equivalent to (1.2).

Conversely, assume that all the coefficients $\alpha_i^{(j)}$ are defined as in (1.2). We denote $c_{00}(\mathbb{N})$ the set of finitely supported sequences $x = (x_k)_{k \in \mathbb{N}}$, which is a dense subspace of $\ell^2(\mathbb{N})$. Let $(e_n)_{n=1}^{+\infty}$ be the canonical orthonormal basis of $\ell^2(\mathbb{N})$ and, for $n \geq 1$, let $E_n = \text{span}(e_1, \ldots, e_n)$ and $P_n$ be the orthogonal projection from $\ell^2(\mathbb{N})$ onto $E_n$.

Let $x \in c_{00}(\mathbb{N})$. For $n$ large enough, we have $P_n x = x$ and, then, $P_n M P_n x = P_n M x$. As we know that $P_n M P_n$ is a contraction, then, for $n$ large enough, we have $\|P_n M x\| \leq \|x\|$, which implies that $\|Mx\| \leq \|x\|$.

Therefore, $M$ defines a contraction on $c_{00}(\mathbb{N})$, which can be extended by density to a contraction $T_M$ on $\ell^2(\mathbb{N})$. The matrix representation of $T_M$ with respect to the canonical orthonormal basis of $\ell^2(\mathbb{N})$ is obviously $M$.

□

## Acknowledgments

The author thanks Catalin Badea for useful comments and encouragements. The author also acknowledge support from the Labex CEMPI (ANR-11-LABX-0007-01).

(A. Renard) UNIV LILLE, CNRS, UMR 8524 - LABORATOIRE PAUL PAINLEVÉ, FRANCE

*E-mail address*: axel.renard@univ-lille.fr